# Deconvolution for an atomic distribution

**Bert van Es**

*Korteweg-de Vries Institute for Mathematics*
*Universiteit van Amsterdam*
*Plantage Muidergracht 24*
*1018 TV Amsterdam*
*The Netherlands*
*e-mail:* vanes@science.uva.nl

**Shota Gugushvili**[*][†]

*Eurandom*
*Technische Universiteit Eindhoven*
*P.O. Box 513*
*5600 MB Eindhoven*
*The Netherlands*
*e-mail:* gugushvili@eurandom.tue.nl

**Peter Spreij**

*Korteweg-de Vries Institute for Mathematics*
*Universiteit van Amsterdam*
*Plantage Muidergracht 24*
*1018 TV Amsterdam*
*The Netherlands*
*e-mail:* spreij@science.uva.nl

**Abstract:** Let $X_1, \ldots, X_n$ be i.i.d. observations, where $X_i = Y_i + \sigma Z_i$ and $Y_i$ and $Z_i$ are independent. Assume that unobservable $Y$'s are distributed as a random variable $UV$, where $U$ and $V$ are independent, $U$ has a Bernoulli distribution with probability of zero equal to $p$ and $V$ has a distribution function $F$ with density $f$. Furthermore, let the random variables $Z_i$ have the standard normal distribution and let $\sigma > 0$. Based on a sample $X_1, \ldots, X_n$, we consider the problem of estimation of the density $f$ and the probability $p$. We propose a kernel type deconvolution estimator for $f$ and derive its asymptotic normality at a fixed point. A consistent estimator for $p$ is given as well. Our results demonstrate that our estimator behaves very much like the kernel type deconvolution estimator in the classical deconvolution problem.

Received September 2007.

---

[*]The corresponding author.
[†]The research of this author was financially supported by the Nederlandse Organisatie voor Wetenschappelijk Onderzoek (NWO). The research was conducted while this author was at Korteweg-de Vries Institute for Mathematics in Amsterdam.

**Contents**

## 1. Introduction

Let $X_1, \ldots, X_n$ be i.i.d. copies of a random variable $X = Y + \sigma Z$, where $X_i = Y_i + \sigma Z_i$, and $Y_i$ and $Z_i$ are independent and have the same distribution as $Y$ and $Z$, respectively. Assume that $Y$'s are unobservable and that $Y = UV$, where $U$ and $V$ are independent, $U$ has a Bernoulli distribution with probability of zero equal to $p$ (we assume that $0 \leq p < 1$) and $V$ has a distribution function $F$ with density $f$. Furthermore, let the random variable $Z$ have a standard normal distribution and let $\sigma$ be a known positive number. The $X$ will then have a density, which we denote by $q$. The distribution of $Y$ is completely determined by $f$ and $p$. Note that the distribution of $Y$ has an atom at zero. Based on a sample $X_1, \ldots, X_n$, we consider the problem of (nonparametric) estimation of the density $f$ and the probability $p$.

Our estimation problem is closely related to the classical deconvolution problem, where the situation is as described above, except that in the classical case $p$ vanishes and $Y_i$ has a continuous distribution with density $f$, which we want to estimate. The $Y_i$'s can for instance be interpreted as measurements of some characteristic of interest, contaminated by noise $\sigma Z_i$. Some works on deconvolution include [3, 4, 6, 7, 9, 10, 11, 13, 14, 16, 19, 20, 21, 22, 23, 28, 30, 32, 35, 38, 39, 42, 43, 45, 46] and [50]. Practical problems related to deconvolution can be found e.g. in [31], which provides a general account of mixture models. The deconvolution problem is also related to empirical Bayes estimation of the prior distribution, see e.g. [2] and [33]. Yet another application field is the nonparametric errors in variables regression, see [24].

Unlike the classical deconvolution problem, in our case $Y$ does not have a density, because the distribution of $Y$ has an atom at zero. Hence our results, apart of the direct applications below, will also provide insight into the robustness of the deconvolution estimator when the assumption of absolute continuity is violated.

One situation where the atomic deconvolution can arise, is the following: one might think of the $X_i$'s as increments $\mathbf{X}_i - \mathbf{X}_{i-1}$ of a stochastic process $\mathbf{X}_t = \mathbf{Y}_t + \sigma \mathbf{Z}_t$, where $\mathbf{Y} = (\mathbf{Y}_t)_{t \geq 0}$ is a compound Poisson process with intensity $\lambda$ and jump size density $\rho$, and $\mathbf{Z} = (\mathbf{Z}_t)_{t \geq 0}$ is a Brownian motion independent of $\mathbf{Y}$. The distribution of $\mathbf{Y}_i - \mathbf{Y}_{i-1}$ then has an atom at zero with probability equal to $e^{-\lambda}$, while $\mathbf{Z}_i - \mathbf{Z}_{i-1}$ has a standard normal distribution. Notice that $\mathbf{X} = (\mathbf{X}_t)_{t\geq 0}$ is a Lévy process, see Example 8.5 in [37]. An exponential of the process $\mathbf{X}$ can be used to model the evolution of a stock price, see [34]. The law of $\mathbf{X}$ can be completely characterised by $f$, $\lambda$ and $\sigma$. Furthermore, estimation of $f$ in the atomic deconvolution context is closely related to estimation of the jump size density of a compound Poisson process $\mathbf{Y}$, which is contaminated by noise coming from a Brownian motion, see [26].

Another practical situation might arise in missing data problems. Suppose for instance that a measurement device is used to measure some quantity of interest and that it has a fixed probability $p$ of failure to detect this quantity, in which case it renders zero. Repetitive measurements can be modelled by random variables $Y_i$ defined as above. Assume that our goal is to estimate the density $f$ and the probability $p$. In practice measurements are often contaminated by an additive measurement error and to account for this, we add the noise $\sigma Z_i$ to our measurements ($\sigma$ quantifies the noise level). If we could directly use the measurements $Y_i$, then the zero measurements could be discarded and we would have observations with density $f$ to base our estimator on. However, due to the additional noise $\sigma Z_i$, the zeroes cannot be distinguished from the nonzero $Y_i$'s. The use of deconvolution techniques is thus unavoidable. The same situation occurs for instance when $Y_i$ are left truncated at zero. In the error-free case, i.e. when $\sigma = 0$, estimation of the mean and variance of a *positive* random variable $V$ was considered in [1]. Our model appears to be more general.

In what follows, we first assume that $p$ is known and construct an estimator for $f$. After this, in the model where $p$ is unknown, we will provide an estimator for $p$ and then propose a plug-in type estimator for $f$. An estimator for $f$ will be constructed via methods similar to those used in the classical deconvolution problem. In particular we will use Fourier inversion and kernel smoothing. Let $\phi_X$, $\phi_Y$ and $\phi_f$ denote the characteristic functions of the random variables $X, Y$ and $V$, respectively. Notice that the characteristic function of $Y$ is given by

$$\phi_Y(t) = p + (1-p)\phi_f(t). \tag{1.1}$$

Furthermore, since

$$\phi_X(t) = \phi_Y(t)e^{-\sigma^2 t^2/2} = (p + (1-p)\phi_f(t))e^{-\sigma^2 t^2/2},$$

the characteristic function of $V$ can be expressed as

$$\phi_f(t) = \frac{\phi_X(t) - pe^{-\sigma^2 t^2/2}}{(1-p)e^{-\sigma^2 t^2/2}}.$$

Assuming that $\phi_f$ is integrable, by Fourier inversion we get

$$f(x) = \frac{1}{2\pi} \int_{-\infty}^{\infty} e^{-itx} \frac{\phi_X(t) - pe^{-\sigma^2 t^2/2}}{(1-p)e^{-\sigma^2 t^2/2}} dt. \tag{1.2}$$

An obvious way to construct an estimator of $f(x)$ from this relation is to estimate the characteristic function $\phi_X(t)$ by its empirical counterpart,

$$\phi_{emp}(t) = \frac{1}{n} \sum_{j=1}^{n} e^{itX_j},$$

see e.g. [25] for a discussion of its applications in statistics, and then obtain the estimator of $f$ by a plug-in device. Alternatively, one can estimate the density $q$ of $X$ by a kernel estimator

$$q_{nh}(x) = \frac{1}{nh}\sum_{j=1}^{n} w\left(\frac{x-X_j}{h}\right),$$

where $w$ denotes a kernel function and $h > 0$ is a bandwidth. Denote by $\phi_w$ the Fourier transform of the kernel $w$. The characteristic function of $q_{nh}$, which is equal to $\phi_{emp}(t)\phi_w(ht)$, will serve as an estimator of $\phi_q$, the characteristic function of $q$. A naive estimator of $f$ can then be obtained by a plug-in device, and would be

$$\frac{1}{2\pi}\int_{-\infty}^{\infty} e^{-itx}\frac{\phi_{emp}(t)\phi_w(ht) - pe^{-\sigma^2 t^2/2}}{(1-p)e^{-\sigma^2 t^2/2}}dt. \tag{1.3}$$

However, this procedure is not always meaningful, because the integrand in (1.3) is not integrable in general. Therefore, instead of (1.3), we define our estimator of $f$ as

$$f_{nh}(x) = \frac{1}{2\pi}\int_{-\infty}^{\infty} e^{-itx}\frac{\phi_{emp}(t) - pe^{-\sigma^2 t^2/2}}{(1-p)e^{-\sigma^2 t^2/2}}\phi_w(ht)dt, \tag{1.4}$$

where the integral is well-defined under the assumption that $\phi_w$ has a compact support on $[-1, 1]$. Notice that

$$f_{nh}(x) = \frac{\hat{f}_{nh}(x)}{1-p} - \frac{p}{1-p}w_h(x), \tag{1.5}$$

where

$$\hat{f}_{nh}(x) = \frac{1}{2\pi}\int_{-\infty}^{\infty} e^{-itx}\frac{\phi_{emp}(t)\phi_w(ht)}{e^{-\sigma^2 t^2/2}}dt \tag{1.6}$$

and $w_h(x) = (1/h)w(x/h)$. Hence $\hat{f}_{nh}$ has the same form as an ordinary deconvolution kernel density estimator based on the sample $X_1,\ldots,X_n$, see e.g. pp. 231–232 in [49].

Under the assumption of integrability of $\phi_f$ and some additional restrictions on $w$, the bias of the estimator (1.4) will asymptotically vanish as $h \to 0$. Indeed,

$$\mathrm{E}[f_{nh}(x)] - f(x) = \frac{1}{2\pi}\int_{-\infty}^{\infty} e^{-itx}\phi_f(t)(\phi_w(ht) - 1)dt. \tag{1.7}$$

The result follows via the dominated convergence theorem, once we know that $\phi_w$ is bounded and $\phi(0) = 1$. Observe that (1.7) coincides with the bias of an ordinary kernel density estimator based on a sample from $f$. In case we know that $f$ belongs to a specific Hölder class, it is possible to derive an order bound for (1.7) in terms of some power of $h$, see Proposition 1.2 in [40]. Further properties of kernel density estimators can be found in [15, 17, 36, 40, 48] and [49].

Estimation of $p$ is not as easy, as it might appear at first sight. Indeed, due to the convolution structure $X = Y + \sigma Z$, the random variable $X$ has a density and the atom in the distribution of $Y$ is not inherited by the distribution of $X$. On the other hand $p$ is identifiable, since

$$\lim_{t \to \infty} \frac{\phi_X(t)}{e^{-\sigma^2 t^2/2}} = \lim_{t \to \infty} \phi_Y(t) = p,$$

because $\phi_f(t) \to 0$ as $t \to \infty$ by the Riemann-Lebesgue theorem. However, this relation cannot be used as a hint for the construction of a meaningful estimator of $p$ because of the oscillating behaviour of $\phi_{emp}(t)$, the obvious estimator of $\phi_X(t)$, as $t \to \infty$.

As an estimator of $p$ we propose

$$p_{ng} = \frac{g}{2} \int_{-1/g}^{1/g} \frac{\phi_{emp}(t) \phi_k(gt)}{e^{-\sigma^2 t^2/2}} dt, \tag{1.8}$$

where the number $g > 0$ denotes a bandwidth and $\phi_k$ denotes the Fourier transform of a kernel $k$. We assume that $\phi_k$ has support $[-1, 1]$. The definition of $p_{ng}$ is motivated by the fact that

$$\lim_{g \to 0} \frac{g}{2} \int_{-1/g}^{1/g} \frac{\phi_X(t)}{e^{-\sigma^2 t^2/2}} dt = \lim_{g \to 0} \frac{g}{2} \int_{-1/g}^{1/g} \phi_Y(t) dt$$

$$= \lim_{g \to 0} \frac{g}{2} \int_{-1/g}^{1/g} (p + (1-p)\phi_f(t)) dt$$

$$= p.$$

Assuming the integrability of $\phi_f$, the last equality follows from

$$\int_{-1/g}^{1/g} |\phi_f(t)| dt \leq \int_{-\infty}^{\infty} |\phi_f(t)| dt < \infty.$$

Finally, let us consider the general case when both $p$ and $f$ are unknown. Plugging in an estimator of $p$ into (1.4) leads to the following definition of an estimator of $f$,

$$f_{nhg}^*(x) = \frac{1}{2\pi} \int_{-\infty}^{\infty} e^{-itx} \frac{\phi_{emp}(t) - \hat{p}_{ng} e^{-\sigma^2 t^2/2}}{(1 - \hat{p}_{ng}) e^{-\sigma^2 t^2/2}} \phi_w(ht) dt, \tag{1.9}$$

where

$$\hat{p}_{ng} = \min(p_{ng}, 1 - \epsilon_n). \tag{1.10}$$

Here $0 < \epsilon_n < 1$ and $\epsilon_n \downarrow 0$ at a suitable rate, which will be specified in Condition 1.5. The truncation of $p_{ng}$ in (1.10) is introduced for technical reasons, see formula (5.18), where we need that the random variable $1 - \hat{p}_{ng}$ is bounded away from zero.

In practice it might also happen that the error variance $\sigma^2$ is unknown and hence has to be estimated. This is a difficult problem in the classical deconvolution density estimation if only observations $X_1, \ldots, X_n$ are available, as the convergence rate for estimation of $\sigma$ is not the usual $\sqrt{n}$ rate, see e.g. [35]. Moreover, the convergence rate of an estimator of $\sigma$ would dominate the asymptotics. If additional measurements are available, then as suggested for instance in [9], $\sigma^2$ can be estimated e.g. via the empirical variance of the difference of replicated observations or by the method of moments via instrumental variables. A recent paper on this subject is [14]. We do not pursue this question any further and assume that $\sigma$ is known.

Concluding this section, we introduce some technical conditions on the density $f$, kernels $w$ and $k$, bandwidths $h$ and $g$ and the sequence $\epsilon_n$. These are needed in the proof of Theorem 2.5, the main theorem of the paper, and subsequent results. Weaker forms of these conditions are sufficient to prove other results from Section 2 and will be given directly in the corresponding statements.

**Condition 1.1.** *There exists a number $\gamma > 0$, such that $u^\gamma \phi_f(u)$ is integrable.*

**Condition 1.2.** *Let $\phi_w$ be bounded, real valued, symmetric and have support $[-1, 1]$. Let $\phi_w(0) = 1$ and let*

$$\phi_w(1-t) = At^\alpha + o(t^\alpha), \qquad as \qquad t \downarrow 0 \tag{1.11}$$

*for some constants $A$ and $\alpha \geq 0$. Moreover, we assume that $\gamma > 1 + 2\alpha$.*

This condition is similar to the one used in [30] and [46] in the classical deconvolution problem. An example of a kernel that satisfies this condition is

$$w(x) = -\frac{4\sqrt{2}(3x\cos x + (-3 + x^2)\sin x)}{\pi x^5}. \tag{1.12}$$

Its Fourier transform is given by

$$\phi_w(t) = (1-t^2)^2 1_{[-1,1]}(t). \tag{1.13}$$

In this case $\alpha = 2$ and $A = 4$. The kernel (1.12) and its Fourier transform are plotted in Figures 1 and 2.

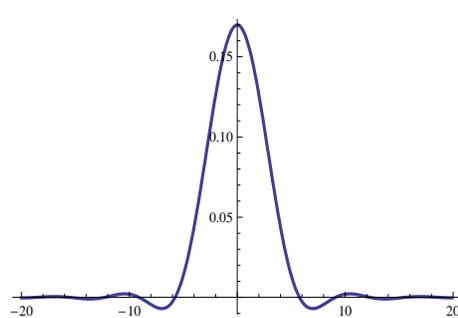

FIG 1. *The kernel (1.13).*

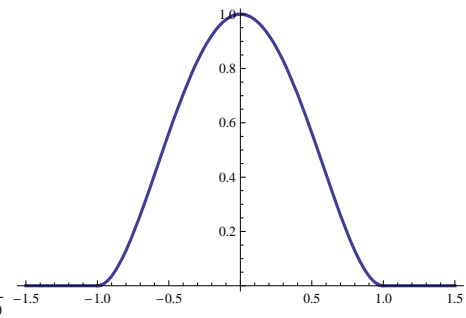

FIG 2. *The Fourier transform of the kernel (1.13).*

**Condition 1.3.** *Let $\phi_k$ be real valued, symmetric and have support $[-1,1]$. Let $\phi_k$ integrate to 2 and let*

$$\phi_k(t) = Bt^\gamma + o(t^\gamma), \tag{1.14}$$
$$\phi_k(1-t) = Ct^\alpha + o(t^\alpha), \tag{1.15}$$

*as $t \downarrow 0$. Here $B$ and $C$ are some constants, and $\gamma$ and $\alpha$ are the same as above.*

An example of such a kernel is given by

$$k(x) = -\frac{2079x(-151200 + 21840x^2 - 730x^4 + 7x^6)\cos x}{\sqrt{2\pi}x^{11}}$$
$$-\frac{693(453600 - 216720x^2 + 13950x^4 - 255x^6 + x^8)\sin x}{\sqrt{2\pi}x^{11}}. \tag{1.16}$$

Its Fourier transform is given by

$$\phi_k(t) = \frac{693}{8}t^6(1-t^2)^2 1_{[-1,1]}(t). \tag{1.17}$$

In this case $B = 693/8$, $\gamma = 6$, $\alpha = 2$ and $C = 693/2$. The kernel (1.16) and its Fourier transform are plotted in Figures 3 and 4. Condition (1.14) is only needed when $p_{ng}$ is plugged into $f^*_{nhg}$, but not if $p_{ng}$ is used as an estimator of $p$.

**Condition 1.4.** *Let the bandwidths $h$ and $g$ depend on $n$, $h = h_n$ and $g = g_n$, and let*

$$h_n = \sigma((1+\eta_n)\log n)^{-1/2},$$
$$g_n = \sigma((1+\delta_n)\log n)^{-1/2},$$

*where $\eta_n$ and $\delta_n$ are such that $\eta_n \downarrow 0, \delta_n \downarrow 0, \eta_n - \delta_n > 0$, and*

$$(\eta_n - \delta_n)\log n \to \infty.$$

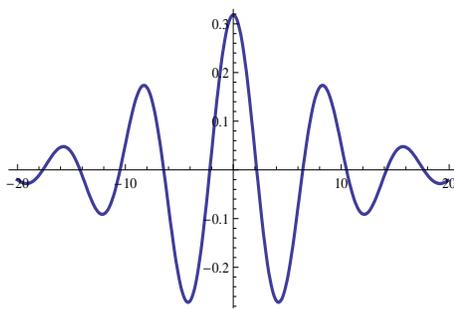
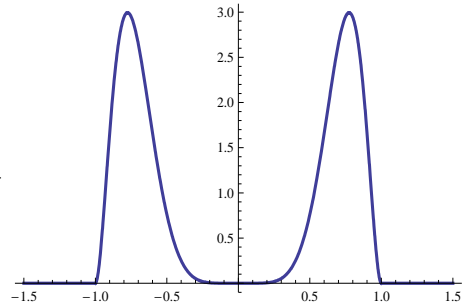

Fig 3. *The kernel (1.17).*

Fig 4. *The Fourier transform of the kernel (1.17).*

*Furthermore, we assume that*

$$-\eta_n \log n + (1+2\alpha) \log \log n \to \infty,$$
$$-\delta_n \log n + (1+2\alpha) \log \log n \to \infty. \tag{1.18}$$

An example of $\eta_n$ and $\delta_n$ in the definition above is

$$\eta_n = 2\frac{\log \log \log n}{\log n}, \quad \delta_n = \frac{\log \log \log n}{\log n}.$$

Conditions on the bandwidths $h_n$ and $g_n$ in Condition 1.4 are not the only possible ones and other restrictions are also possible. However the logarithmic decay of $h_n$ and $g_n$ is unavoidable. Following the default convention in kernel density estimation and to keep the notation compact, we will suppress the index $n$ when writing $h_n$ and $g_n$ and will write $h$ and $g$ instead, since no ambiguity will arise.

**Condition 1.5.** *Let $\epsilon_n \downarrow 0$ be such that*

$$\frac{\log \epsilon_n}{(\eta_n - \delta_n) \log n} \to 0.$$

An example of such $\epsilon_n$ for $\eta_n$ and $\delta_n$ given above is $(\log \log \log n)^{-1}$.

The remainder of the paper is organised as follows: in Section 2 we derive the theorem establishing the asymptotic normality of $f_{nh}(x)$, the fact that the estimator $p_{ng}$ is weakly consistent, and finally that the estimator $f^*_{nhg}(x)$ is asymptotically normal. Section 3 contains simulation examples. Section 4 discusses a method for implementation of the estimator in practice. All the proofs are collected in Section 5.

## 2. Main results

We will first study the estimation of $f$ when $p$ is known, and then proceed to the general case with unknown $p$. The reason for this is twofold. Firstly, it is interesting to compare the behaviour of the estimator of $f$ under the assumption of known and unknown $p$, and secondly, the proofs of the results for the latter case rely heavily on the proofs for the former case.

The first result in this section deals with the nonrobustness of the estimator $\hat{f}_{nh}$. In ordinary kernel deconvolution, when it is assumed that $Y$ is absolutely continuous, the estimator for its density is defined as

$$\hat{f}_{nh}(x) = \frac{1}{2\pi} \int_{-\infty}^{\infty} e^{-itx} \frac{\phi_{emp}(t)\phi_w(ht)}{e^{-\sigma^2 t^2/2}} dt. \tag{2.1}$$

Now suppose that the assumption of absolute continuity of $Y$ is violated. What will happen, if we still use the estimator $\hat{f}_{nh}(x)$? The following result addresses this question.

**Theorem 2.1.** *Let $\hat{f}_{nh}(x)$ be defined as in (2.1). Assume that $\phi_w$ is bounded and has a compact support on $[-1,1]$. Then*

$$\mathrm{E}\,[\hat{f}_{nh}(x)] = p w_h(x) + (1-p) f * w_h(x), \tag{2.2}$$

*where $w_h(\cdot) = (1/h) w(\cdot/h)$, and $*$ denotes convolution.*

From this theorem it follows that $\mathrm{E}\,[\hat{f}_{nh}(0)]$ diverges to infinity as $h \to 0$, because so does $h^{-1} w(0)$, if $w(0) \neq 0$ (the latter is the case for the majority of conventional kernels). In practice this will also result in an equally undesirable behaviour of $\mathrm{E}\,[\hat{f}_{nh}(x)]$ in the neighbourhood of zero. When $x \neq 0$, with a proper selection of a kernel $w$, one can achieve that the first term in (2.2) asymptotically vanishes as $h \to 0$. Indeed, it is sufficient to assume that $w$ is such that $\lim_{u \to \pm \infty} u w(u) = 0$. The second term in (2.2) will converge to $(1-p) f(x)$ as $h \to 0$, provided that $\phi_f$ is integrable, $\phi_w$ is bounded and $\phi_w(0) = 1$. These facts address the issue of the nonrobustness of $\hat{f}_{nh}$: under a misspecified model, i.e. under the assumption that the distribution of $Y$ is absolutely continuous, while in fact it has an atom at zero, the classical deconvolution estimator will exhibit unsatisfactory behaviour near zero. This will happen despite the fact that $\hat{f}_{nh}(x)$ will be asymptotically normal when centred at its expectation and suitably normalised, see Corollary 5.1 in Section 5. The asymptotic normality follows from Lemmas 5.2 and 5.3 of Section 5, where only absolute continuity of the distribution of $X$ is required.

Our next goal is to establish the asymptotic normality of the estimator $f_{nh}(x)$. We formulate the corresponding theorem below.

**Theorem 2.2.** *Assume that $\phi_f$ is integrable. Let $\mathrm{E}\,[X^2] < \infty$, and suppose that Condition 1.2 holds. Let $f_{nh}$ be defined as in (1.4). Then, as $n \to \infty$ and $h \to 0$,*

$$\frac{\sqrt{n}}{h^{1+2\alpha} e^{\sigma^2/(2h^2)}} (f_{nh}(x) - \mathrm{E}\,[f_{nh}(x)])$$
$$\xrightarrow{\mathcal{D}} N\left(0, \frac{A^2}{2\pi^2(1-p)^2} \left(\frac{1}{\sigma^2}\right)^{2+2\alpha} (\Gamma(\alpha+1))^2\right), \tag{2.3}$$

*where $\Gamma$ denotes the gamma function, $\Gamma(t) = \int_0^\infty v^{t-1} e^{-v} dv$.*

Note that Theorem 2.2 establishes asymptotic normality of $f_{nh}$ under an atomic distribution, which constitutes a generalisation of a result in [46] (see also [45]) for the case of the classical deconvolution problem. The generalisation is possible, because the proof uses only the continuity of the density of $X$, which is still true when $Y$ has a distribution with an atom. Furthermore, notice that in order to get a consistent estimator, from this theorem it follows that $\sqrt{n} h^{-1-2\alpha} e^{-\sigma^2/(2h^2)}$ has to diverge to infinity. Therefore the bandwidth $h$ has to be at least of order $(\log n)^{-1/2}$, as it is actually stated in Condition 1.4. In practice this implies that the bandwidth $h$ has to be selected fairly large, even for large sample sizes. This is the case for the classical deconvolution problem as well in the case of a supersmooth error distribution, cf. [46].

Observe that the asymptotic variance in (2.3) does not depend on the target density $f$ nor on the point $x$. This phenomenon is quite peculiar, but is already known in the classical deconvolution kernel density estimation, see for instance equation (6) in [3]. There, provided that $h$ is small enough, the asymptotic variance of the deconvolution kernel density estimator (or, strictly speaking an upper bound for it) also does not depend neither on the target density $f$, nor on the point $x$. In this respect see also [46]. Such results do not contradict the asymptotic normality result in [21], see Theorem 2.2 in that paper, as there the asymptotic variance of the deconvolution kernel density estimator is not evaluated.

Now we state a theorem concerning the consistency of $p_{ng}$, the estimator of $p$.

**Theorem 2.3.** *Assume that $\phi_f$ is integrable, let $\mathrm{E}\left[X^2\right] < \infty$, and let the kernel $k$ have a Fourier transform $\phi_k$ that is bounded by one and integrates to two. Let $p_{ng}$ be defined as in (1.8). If $g$ is such that $g^{4+4\alpha}e^{\sigma^2/g^2}n^{-1} \to 0$, then $p_{ng}$ is a consistent estimator of $p$, i.e.*

$$\mathrm{P}(|p_{ng} - p| > \epsilon) \to 0$$

*as $n \to 0$ and $g \to 0$. Here $\epsilon$ is an arbitrary positive number. Furthermore, under Condition 1.1 and 1.3*

$$\mathrm{E}\left[(p_{ng} - p)^2\right] = O\left(g^{2+2\gamma} + \frac{g^{4+4\alpha}e^{\sigma^2/g^2}}{n}\right).$$

One can also show that $p_{ng}$ is asymptotically normal, when centred and suitably normalised. We formulate the corresponding theorem below.

**Theorem 2.4.** *Assume that the conditions of Theorem 2.3 hold. Let $p_{ng}$ be defined as in (1.8) and let (1.15) hold. Then*

$$\frac{\sqrt{n}}{g^{2+2\alpha}e^{\sigma^2/(2g^2)}}(p_{ng} - \mathrm{E}\left[p_{ng}\right]) \xrightarrow{\mathcal{D}} N\left(0, \frac{C^2(\Gamma(1+\alpha))^2}{2}\left(\frac{1}{\sigma^2}\right)^{2+2\alpha}\right)$$

*as $n \to \infty$ and $g \to 0$.*

Finally, we consider the case when both $p$ and $f$ are unknown. We state the main theorem of the paper.

**Theorem 2.5.** *Let $f^*_{nhg}(x)$ be defined by (1.9), $\mathrm{E}\left[X^2\right] < \infty$ and let Conditions 1.1–1.5 hold. Then, as $n \to \infty$, we have*

$$\frac{\sqrt{n}}{h^{1+2\alpha}e^{\sigma^2/(2h^2)}}(f^*_{nhg}(x) - \mathrm{E}\left[f^*_{nhg}(x)\right])$$
$$\xrightarrow{\mathcal{D}} N\left(0, \frac{A^2}{2\pi^2(1-p)^2}\left(\frac{1}{\sigma^2}\right)^{2+2\alpha}(\Gamma(\alpha+1))^2\right).$$

Notice that the asymptotic variance is the same as in (2.3), which justifies the plug-in approach to the construction of an estimator of $f$, when $p$ is unknown.

A natural question to consider is what happens when we centre $f^*_{nhg}(x)$ not at its expectation, but at $f(x)$. This has practical importance as well, e.g. for the construction of (asymptotic) confidence intervals. Writing

$$\frac{\sqrt{n}}{h^{1+2\alpha}e^{\sigma^2/(2h^2)}}(f^*_{nhg}(x) - f(x)) = \frac{\sqrt{n}}{h^{1+2\alpha}e^{\sigma^2/(2h^2)}}(f^*_{nhg}(x) - \mathrm{E}\,[f^*_{nhg}(x)])$$
$$+ \frac{\sqrt{n}}{h^{1+2\alpha}e^{\sigma^2/(2h^2)}}(\mathrm{E}\,[f^*_{nhg}(x)] - f(x)),$$

we see, that we have to study the second term here, i.e. to compare the behaviour of the bias of $f^*_{nhg}(x)$ to the normalising factor $\sqrt{n}h^{-(1+2\alpha)}e^{-\sigma^2/(2h^2)}$. We will study the bias of $f^*_{nhg}(x)$ in two steps: first we will show that it asymptotically vanishes, which itself is of independent interest. After this we will provide conditions under which it asymptotically vanishes when multiplied by $\sqrt{n}h^{-(1+2\alpha)}e^{-\sigma^2/(2h^2)}$. Recall the definition of a Hölder class of functions $\mathcal{H}(\beta, L)$.

**Definition 2.1.** *A function $f$ is said to belong to the Hölder class $\mathcal{H}(\beta, L)$, if its derivatives up to order $l = [\beta]$ exist and verify the condition*

$$|f^{(l)}(x+t) - f^{(l)}(x)| \leq L|t|^{\beta-l}$$

*for all $x, t \in \mathbb{R}$.*

Such a smoothness condition on a target density $f$ is standard in kernel density estimation, see e.g. p. 5 of [40]. Often one assumes that $\beta = 2$. If $l = 0$, then set $f^{(l)} = f$. We also need the definition of a kernel of order $l$. In particular, we will use the version given in Definition 1.3 of [40].

**Definition 2.2.** *A kernel $w$ is said to be a kernel of order $l$ for $l \geq 1$, if the functions $x \mapsto x^j w(x)$ are integrable for $j = 0, \ldots, l$ and if*

$$\int_{-\infty}^{\infty} w(x)dx = 1, \quad \int_{-\infty}^{\infty} x^j w(x)dx = 0 \quad \text{for } j = 1, \ldots, l.$$

**Theorem 2.6.** *Let $f^*_{nhg}(x)$ be defined by (1.9) and assume conditions of Theorem 2.5. Then, as $n \to \infty$, we have*

$$\mathrm{E}\,[f^*_{nhg}(x)] - f(x) \to 0.$$

*If additionally $f \in \mathcal{H}(\beta, L)$, $w$ is a kernel of order $l = [\beta]$ and $\beta > 1 + 2\alpha$, then*

$$\frac{\sqrt{n}}{h^{1+2\alpha}e^{\sigma^2/(2h^2)}}(\mathrm{E}\,[f^*_{nhg}(x)] - f(x)) \to 0$$

*as $n \to \infty$.*

Combination of this theorem with Theorem 2.5 leads to the following result.

**Theorem 2.7.** *Assume that the conditions of Theorem 2.6 hold. Then, as $n \to \infty$, we have*

$$\frac{\sqrt{n}}{h^{1+2\alpha}e^{\sigma^2/(2h^2)}}(f^*_{nhg}(x) - f(x)) \xrightarrow{\mathcal{D}} N\left(0, \frac{A^2}{2\pi^2(1-p)^2}\left(\frac{1}{\sigma^2}\right)^{2+2\alpha}(\Gamma(\alpha+1))^2\right).$$

One should keep in mind that these results deal only with asymptotics. In the next section we will study several simulation examples, which will provide some insight into the finite sample properties of the estimator.

## 3. Simulation examples

In this section we consider a number of simulation examples. We do not pretend to provide an exhaustive simulation study, rather an illustration, which requires further verification.

Assume that $\sigma = 1$, $p = 0.1$ and that $f$ is normal with mean 3 and variance 9. This results in a nontrivial deconvolution problem, because the ratio of 'noise' compared to 'signal' is reasonably high: $\text{NSR} = \text{Var}[\sigma Z]/\text{Var}[Y]100\% \approx 11\%$. We have simulated a sample of size $n = 1000$. As kernels $w$ and $k$ we selected kernels (1.12) and (1.16), respectively. The bandwidths $h = 0.58$ and $g = 0.5$ were selected by hand. A possible method of computing the estimate is given in Section 4. The estimator $p_{ng}$ produced a value equal to 0.11. The estimate of $f$ (bold dotted line), resulting from the procedure described above, together with the target density $f$ (dashed line) is plotted in Figure 5. For comparison purposes, we have also plotted the estimate $f_{nh}(x)$ (it can be obtained using (1.5) and the true value of the parameter $p$), see Figure 6. As can be seen from the comparison of these two figures, the estimates $f^*_{nhg}$ and $f_{nh}$ look rather similar.

As the second example we consider the case when $f$ is a gamma density with parameters $\alpha = 8$ and $\beta = 1$, i.e.

$$f(x) = \frac{x^7 e^{-x}}{\Gamma(8)} 1_{[x>0]}, \tag{3.1}$$

and $p = 0.25$. We simulated a sample of size $n = 1000$. The kernels were chosen as above and the bandwidths $g = 0.6$ and $h = 0.6$ were selected by hand. The estimate $p_{ng}$ took a value approximately equal to 0.23. The resulting estimate $f^*_{nhg}$ is plotted in Figure 7. As above we also plotted the estimate $\hat{f}_{nh}$, see Figure 8 (notice that the estimate takes on negative values in the neighbourhood of zero). Again both figures look similar.

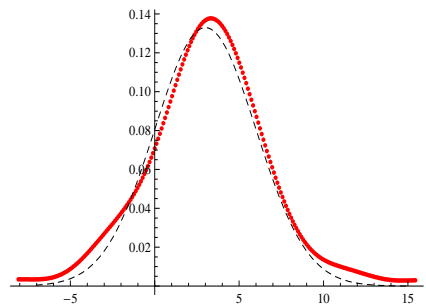 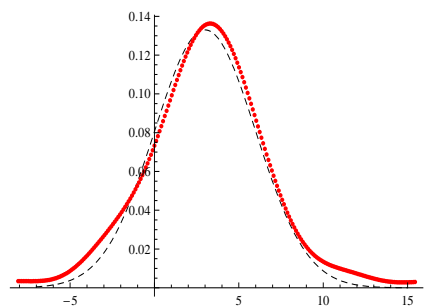

Fig 5. *The normal density $f$ (dashed line) and the estimate $f^*_{nhg}$ (solid line). The sample size $n = 1000$.*

Fig 6. *The normal density $f$ (dashed line) and the estimate $f_{nh}$ (solid line). The sample size $n = 1000$.*

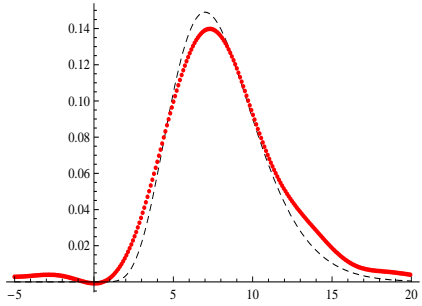

Fig 7. *The gamma density (dashed line) and the estimate $f^*_{nhg}$ (solid line). The sample size $n = 1000$.*

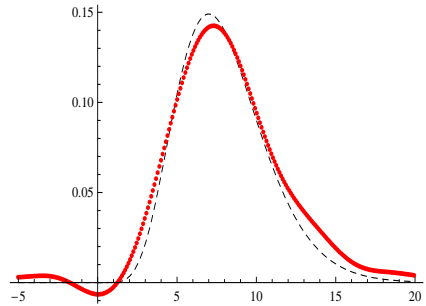

Fig 8. *The gamma density (dashed line) and the estimate $f_{nh}$ (solid line). The sample size $n = 1000$.*

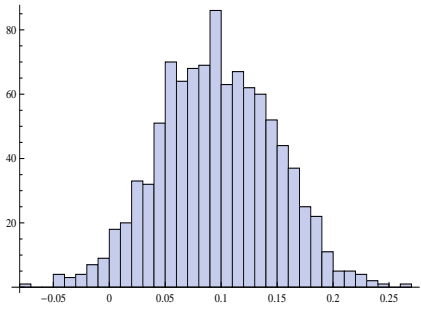

Fig 9. *The histogram of estimates of p for $g = 0.5$ and the sample size $n = 1000$.*

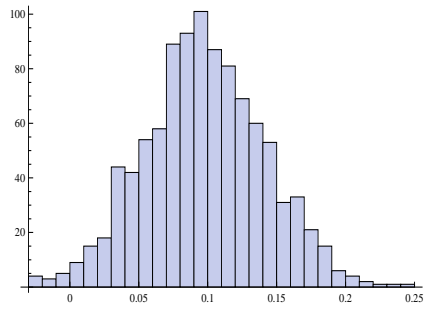

Fig 10. *The histogram of estimates of p for $g = 0.55$ and the sample size $n = 1000$.*

Examination of these figures leads us to two questions: how well does $p_{ng}$ estimate $p$ for moderate sample samples? How sensitive is $f^*_{nhg}$ to under- or overestimation of $p$? To get at least a partial answer to the first question, we considered the same model as in our first example in this section (i.e. deconvolution of the normal density) and repeatedly, i.e. 1000 times, estimated $p$ for the bandwidth $g = 0.5$ and the sample size $n = 1000$ for each simulation run. Then the same procedure was repeated for the bandwidths $g = 0.55$, 0.6 and 0.65. The resulting histograms are plotted in Figures 9–12. They look quite satisfactory. The sample means and sample standard deviations (SD) of estimates of $p$ for different choices of bandwidth $g$ together with the theoretical standard deviations are summarised in Table 1. One notices that the sample means in Table 1 are close to the true value 0.1 of the parameter $p$. The theoretical standard deviations in the same table were computed using Theorem 2.4, which predicts

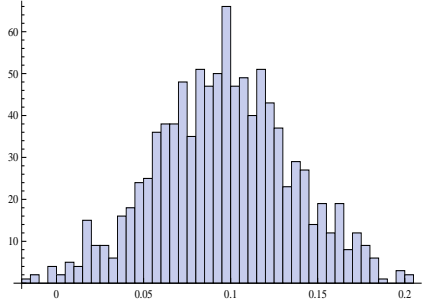
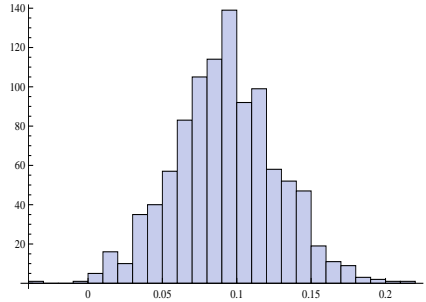

Fig 11. *The histogram of estimates of p for $g = 0.6$ and the sample size $n = 1000$.*

Fig 12. *The histogram of estimates of p for $g = 0.65$ and the sample size $n = 1000$.*

Table 1
*Sample and theoretical means and standard deviations (SD) of estimates of p for different choices of bandwidth g. The sample size $n = 1000$*

| Bandwidth | 0.5 | 0.55 | 0.6 | 0.65 |
|---|---|---|---|---|
| Sample mean | 0.0963 | 0.0975 | 0.0960 | 0.0927 |
| Sample SD | 0.0516 | 0.0436 | 0.0388 | 0.0349 |
| Asymptotic SD | 1.7891 | 2.2399 | 2.8994 | 3.8164 |
| Theoretical SD | 0.0700 | 0.0593 | 0.0487 | 0.0432 |

that they should be equal to (recall, that in our case $\alpha = 2$)

$$\frac{g^6 e^{1/(2g^2)}}{\sqrt{n}} C\sqrt{2}.$$

From Table 1 one sees that there is a large discrepancy between the sample standard deviations and the standard deviations predicted by the theory. The explanation of this discrepancy lies in the fact that the proof of the asymptotic normality of $p_{ng}$ heavily relies on the asymptotic equivalence

$$\int_0^1 \phi_k(s) e^{\sigma^2 s^2/(2g^2)} ds \sim C\Gamma(1+\alpha)\left(\frac{1}{\sigma^2}\right)^{1+\alpha} g^{2(1+\alpha)} e^{\sigma^2/(2g^2)}, \qquad (3.2)$$

see Lemma 5.1 and the proof of Lemma 5.2 in Section 5 below. However, by direct evaluation of the integral on the left-hand side of (3.2) for different values of $g$, it can be seen that this relation does not provide an accurate approximation in those cases where the bandwidth is relatively large, as it actually is in our case. It then follows that the asymptotic standard deviation will not provide a good approximation of the sample standard deviation unless the bandwidth is very small. This in turn implies that the corresponding sample size must be extremely large. We can correct for this poor approximation of the integral in (3.2) by using the integral itself as a normalising factor instead of the right-hand side of (3.2). The results of this correction are represented in the last line of Table 1. As it can be seen, the theoretical standard deviation and the sample standard deviation are much closer to each other. Since the kernel $k$ was selected more or less arbitrarily, one is tempted to believe that an inaccurate approximation in (3.2) might be due to the kernel. This might be the case, however to a certain degree this seems to be characteristic of all popular kernels employed in kernel deconvolution. Consider for instance the kernel

$$w(x) = \frac{48x(x^2-15)\cos x - 144(2x^2-5)\sin x}{\pi x^7}. \qquad (3.3)$$

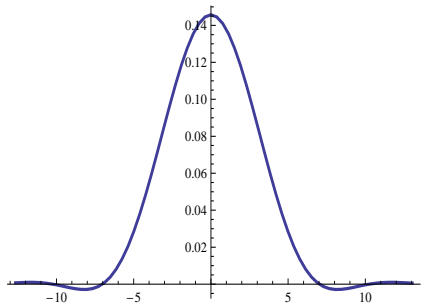

FIG 13. *The kernel* (3.3).

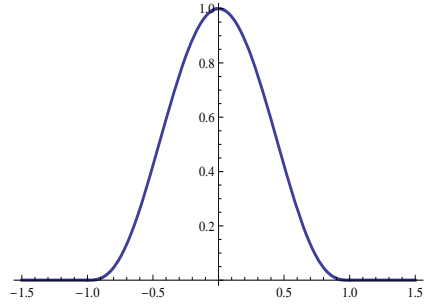

FIG 14. *The Fourier transform of the kernel* (3.3).

Its Fourier transform is given by

$$\phi_w(t) = (1-t^2)^3 1_{[|t|<1]}.$$

The kernel $w$ and its Fourier transform are plotted in Figures 13 and 14, respectively. This kernel was used for simulations in [22] and [47] and it was shown in [13] that it performs well in a deconvolution setting. Notice that this kernel cannot be used to estimate $p$ if we want to plug in the resulting estimator $p_{ng}$ into $f^*_{nhg}$. However, this kernel satisfies Condition 1.2 and can be used to estimate $f$. Nevertheless, the ratio of the left and right hand sides in (3.2) for $h = 0.5$ is equal to 0.4299, which is still far from 1. This issue is further discussed in [42]. Another issue here is that often the error variance $\sigma^2$ is quite small and it is sensible to treat $\sigma$ as depending on the sample size $n$ (with $\sigma \to 0$ as $n \to \infty$), see [9]. However, this is a different model and this question is not addressed here. Notice also that a perfect match between the sample standard deviation and the theoretical standard deviation is impossible to obtain, because we neglect a remainder term when computing the latter. How large the contribution of the remainder term can be in general requires a separate simulation study.

We also considered the case when the error term variance and the sample size are smaller (the target density $f$ was again the standard normal density, while $p$ was set to be 0.1). In particular, we took $\sigma = 0.3$ and $n = 500$. The corresponding histograms are given in Figures 15–18, while the sample and theoretical characteristics for four different choices of the bandwidth $g = 0.5, 0.55, 0.6$ and $0.65$ are summarised in Table 2. Notice a particularly bad match between the asymptotic standard deviation and its empirical counterpart. Other conclusions are similar to those in the previous example.

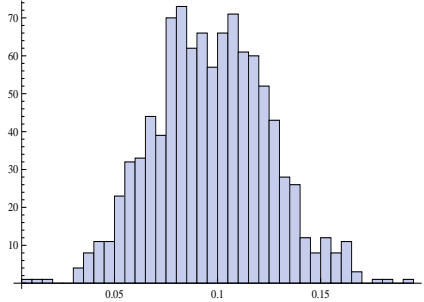

Fig 15. *The histogram of estimates of p for $g = 0.45$ and the sample size $n = 500$.*

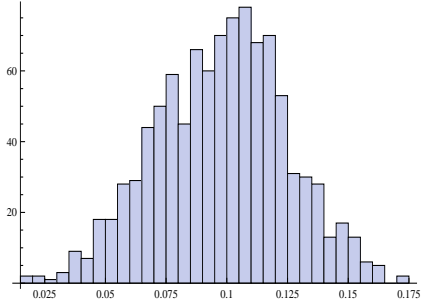

Fig 16. *The histogram of estimates of p for $g = 0.5$ and the sample size $n = 500$.*

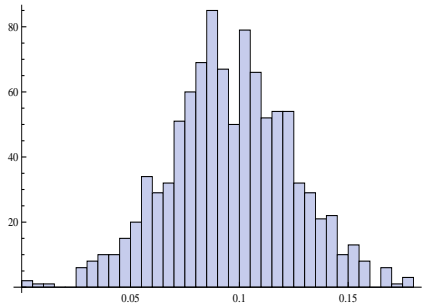

Fig 17. *The histogram of estimates of p for $g = 0.6$ and the sample size $n = 500$.*

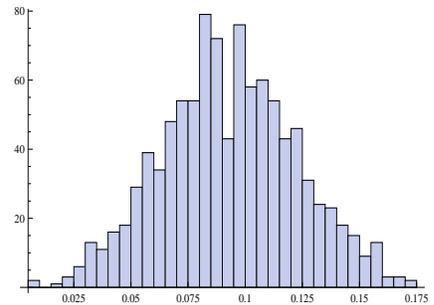

Fig 18. *The histogram of estimates of p for $g = 0.65$ and the sample size $n = 500$.*

To test the robustness of the estimator $f^*_{nhg}$ with respect to the estimated value of $p$, we again turned to the model that was considered in the first example of this section. Instead of $\hat{p}_{ng}$ three different values $\hat{p} = 0.05$, $\hat{p} = 0.1$ and $\hat{p} = 0.15$ were plugged in into (4.2). The resulting estimates $f^*_{nhg}$ are plotted in Figure 19 (the true density is represented by the dashed line). As one can see from Figure 19, under- or overestimation of $p$ in the given range does not have a significant impact on the resulting estimate $\hat{f}_{nh}$ (of course one should keep in mind that $p$ is relatively small in this case). On the other hand, if the value of $\hat{p}$ were larger, e.g. if $\hat{p} = 0.2$, that would have a noticeable effect, e.g. it could have suggested bimodality in the case where the density is actually unimodal, see Figure 20 on the facing page. At the same time the simulated examples concerning the estimates $p_{ng}$ that we considered above seem to suggest that such instances of unsatisfactory estimates of $p$ are not too frequent, because most of the observed values of $p_{ng}$ are concentrated in the interval $[0.05, 0.15]$. We also considered the case when $f = 0.5\phi_{-2,1} + 0.5\phi_{2,1}$, where $\phi_{x,y}$ denotes the normal density with mean $x$ and variance $y$. Hence in this case $f$ is a mixture of two normal densities and it is also bimodal. The match is visually slightly worse for $\hat{p} = 0.2$, but it is still acceptable.

TABLE 2
*Sample and theoretical means and standard deviations (SD) of estimates of p for different choices of bandwidth g. The sample size $n = 500$*

| Bandwidth | 0.45 | 0.5 | 0.6 | 0.65 |
|---|---|---|---|---|
| Sample mean | 0.0972 | 0.0977 | 0.0959 | 0.0930 |
| Sample SD | 0.0277 | 0.0269 | 0.0283 | 0.0295 |
| Asymptotic SD | 311.7 | 562.3 | 2247 | 2521 |
| Theoretical SD | 0.0357 | 0.0349 | 0.0338 | 0.0335 |

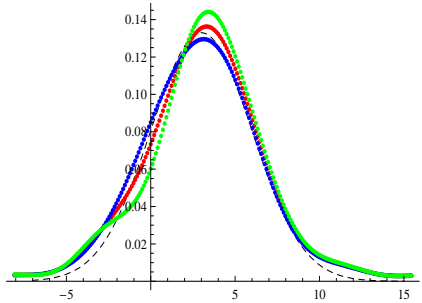
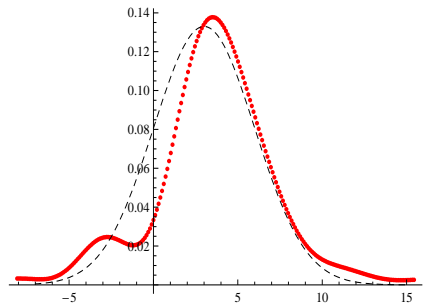

FIG 19. *The normal density $f$ and estimates $f^*_{nhg}$ evaluated for $\hat{p} = 0.05$, $\hat{p} = 0.1$, $\hat{p} = 0.15$ and the sample size $n = 1000$.*

FIG 20. *The normal density $f$ and estimate $f^*_{nhg}$ evaluated for $\hat{p} = 0.2$ and the sample size $n = 1000$.*

The simulation examples that we considered in this section suggest that, despite the slow (logarithmic) rate of convergence, the estimator $f^*_{nhg}$ works in practice (given that $p$ is estimated accurately). This is somewhat comparable to the classical deconvolution problem, where by finite sample calculations it was shown in [47] that for lower levels of noise, the kernel estimators perform well for reasonable sample sizes, in spite of slow rates of convergence for the supersmooth deconvolution problem, obtained e.g. in [21] and [22]. However, Condition 1.4 tells us, that the bandwidths $h$ and $g$ have to be of order $(\log n)^{-1/2}$. In practice this implies that to obtain reasonable estimates, the bandwidths have to be selected fairly large, even for large samples.

One more practical issue concerning the implementation of the estimator $f^*_{nhg}$ (or $p_{ng}$) is the method of bandwidth selection, which is not addressed in this paper. We expect that techniques similar to those used in the classical deconvolution problem will produce comparable results in our problem. This requires a separate investigation of the behaviour of the mean integrated square error of $f^*_{nhg}$. In the case of the classical deconvolution problem papers that consider the issue of data-dependent bandwidth selection are [10, 11, 18, 28] and [39]. Yet another issue is the choice of kernels $w$ and $k$. For the case of the classical deconvolution problem we refer to [13]. In general in kernel density estimation it is thought that the choice of a kernel is of less importance for the performance of an estimator than the choice of the bandwidth, see e.g. p. 31 in [48], or p. 132 in [49].

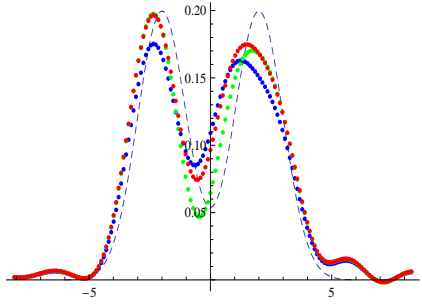
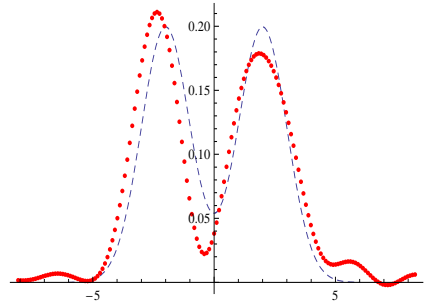

FIG 21. *The mixture of normal densities $f$ and estimates $f^*_{nhg}$ evaluated for $\hat{p} = 0.05$, $\hat{p} = 0.1$, $\hat{p} = 0.15$. and the sample size $n = 1000$.*

FIG 22. *The mixture of normal densities $f$ and estimate $f^*_{nhg}$ evaluated for $\hat{p} = 0.2$ and the sample size $n = 1000$.*

## 4. Computational issues

To compute the estimator $f^*_{nhg}$ in Section 3, a method similar to the one used in [44] (in turn motivated by [5]) can be employed. Namely, notice that

$$\hat{f}_{nh}(x) = \hat{f}^{(1)}_{nh}(x) + \hat{f}^{(2)}_{nh}(x),$$

where

$$\hat{f}^{(1)}_{nh}(x) = \frac{1}{2\pi} \int_0^\infty e^{-itx} \frac{\phi_{emp}(t)\phi_w(ht)}{e^{-\sigma^2 t^2/2}} dt,$$

$$\hat{f}^{(2)}_{nh}(x) = \frac{1}{2\pi} \int_0^\infty e^{itx} \frac{\phi_{emp}(-t)\phi_w(ht)}{e^{-\sigma^2 t^2/2}} dt.$$

Using the trapezoid rule and setting $v_j = \eta(j-1)$, we have

$$\hat{f}^{(1)}_{nh}(x) \approx \frac{1}{2\pi} \sum_{j=1}^{N} e^{-iv_j x} \psi(v_j) \eta, \qquad (4.1)$$

where $N$ is some power of 2 and

$$\psi(v_j) = \frac{\phi_{emp}(v_j)\phi_w(hv_j)}{e^{-\sigma^2 v_j^2/2}}.$$

The Fast Fourier Transform is used to compute values of $\hat{f}^{(1)}_{nh}$ at $N$ different points (concerning the application of the Fast Fourier Transform in kernel deconvolution see [12]). We employ a regular spacing size $\delta$, so that the values of $x$ are

$$x_u = -\frac{N\delta}{2} + \delta(u-1),$$

where $u = 1, \ldots, N$. Therefore, we obtain

$$\hat{f}^{(1)}_{nh}(x_u) \approx \frac{1}{2\pi} \sum_{j=1}^{N} e^{-i\delta\eta(j-1)(u-1)} e^{iv_j \frac{N\delta}{2}} \psi(v_j) \eta.$$

In order to apply the Fast Fourier Transform, note that we must take

$$\delta\eta = \frac{2\pi}{N}.$$

It follows that a small $\eta$, which is needed to achieve greater accuracy in integration, will result in values of $x$ which are relatively far from each other. Therefore, to improve the integration precision, we will apply Simpson's rule, i.e.

$$\hat{f}_{nh}^{(1)}(x_u) \approx \frac{1}{2\pi} \sum_{j=1}^{N} e^{-i\frac{2\pi}{N}(j-1)(u-1)} e^{iv_j \frac{N\delta}{2}} \psi(v_j) \frac{\eta}{3}(3 + (-1)^j - \delta_{j-1}),$$

where $\delta_j$ denotes the Kronecker symbol (recall, that $\delta_j$ is 1, if $j = 0$ and is 0 otherwise). The same reasoning can be applied to $\hat{f}_{nh}^{(2)}(x)$. The estimate $f_{nhg}^*$ can then be computed by noticing that

$$f_{nhg}^*(x) = \frac{\hat{f}_{nh}(x)}{1-\hat{p}_{ng}} - \frac{\hat{p}_{ng}}{1-\hat{p}_{ng}} w_h(x). \tag{4.2}$$

One should keep in mind that even though $w_h$ can be evaluated directly, it is preferable to use the Fast Fourier Transform for its computation, thus avoiding possible numerical issues, see [12]. Also notice that the direct computation of $\phi_{emp}$ is rather time-demanding for large samples. One way to avoid this problem is to use WARPing, cf. [27]. However, for the purposes of the present study, we restricted ourselves to the direct evaluation of $\phi_{emp}$.

## 5. Proofs

*Proof of Theorem 2.1.* The proof is elementary and is based on the definition of $\hat{f}_{nh}(x)$. By Fubini's theorem we have

$$\mathrm{E}\left[\frac{1}{2\pi}\int_{-\infty}^{\infty} e^{-itx}\frac{\phi_{emp}(t)\phi_w(ht)}{e^{-\sigma^2 t^2/2}}dt\right] = \frac{1}{2\pi}\int_{-\infty}^{\infty} e^{-itx}\phi_Y(t)\phi_w(ht)dt.$$

Recalling that $\phi_Y(t) = p + (1-p)\phi_f(t)$, we obtain

$$\mathrm{E}\left[\hat{f}_{nh}(x)\right] = pw_h(x) + (1-p)f * w_h(x). \tag{5.1}$$

Here we used the facts that

$$\frac{1}{2\pi}\int_{-\infty}^{\infty} e^{-itx}\phi_w(ht)dt = w_h(x),$$

$$\frac{1}{2\pi}\int_{-\infty}^{\infty} e^{-itx}\phi_f(t)\phi_w(ht)dt = f * w_h(x).$$

This concludes the proof. □

The proof of Theorem 2.2 is based on the following three lemmas, all of which are reformulations of results from [46].

**Lemma 5.1.** *Assume Condition 1.2. For $h \to 0$ and $\delta \geq 0$ fixed we have*

$$\int_0^1 (1-s)^\delta \phi_w(s) e^{\sigma^2 s^2/(2h^2)} ds \sim A\Gamma(1+\alpha+\delta) \left(\frac{1}{\sigma^2}\right)^{1+\alpha+\delta} h^{2(1+\alpha+\delta)} e^{\sigma^2/(2h^2)}. \tag{5.2}$$

*Proof.* We follow the same line of thought as in [46]. Using the substitution $s = 1 - h^2 v$ and the dominated convergence theorem in the one but last step, we get

$$\int_0^1 (1-s)^\delta \phi_w(s) e^{\sigma^2 s^2/(2h^2)} ds$$
$$= h^2 \int_0^{1/h^2} (h^2 v)^\delta \phi_w(1-h^2v) e^{\sigma^2(1-h^2v)^2/(2h^2)} dv$$
$$= h^2 \int_0^{1/h^2} \frac{\phi_w(1-h^2v)}{(h^2v)^\alpha} (h^2v)^{\alpha+\delta} e^{\sigma^2(1-h^2v)^2/(2h^2)} dv$$
$$= h^{2+2\alpha+2\delta} e^{\sigma^2/(2h^2)} \int_0^{1/h^2} \frac{\phi_w(1-h^2v)}{(h^2v)^\alpha} v^{\alpha+\delta} e^{-\sigma^2 v + \sigma^2 h^2 v^2/2} dv$$
$$\sim h^{2+2\alpha+2\delta} e^{\sigma^2/(2h^2)} A \int_0^\infty v^{\alpha+\delta} e^{-\sigma^2 v} dv$$
$$= h^{2+2\alpha+2\delta} e^{\sigma^2/(2h^2)} \left(\frac{1}{\sigma^2}\right)^{1+\alpha+\delta} A\Gamma(\alpha+\delta+1).$$

The lemma is proved. □

**Lemma 5.2.** *Assume Condition 1.2 and let $\mathrm{E}\left[X^2\right] < \infty$. Furthermore, let $\hat{f}_{nh}$ be defined by (1.6). Then as $n \to \infty$ and $h \to 0$,*

$$\frac{\sqrt{n}}{h^{1+2\alpha} e^{\sigma^2/(2h^2)}} (\hat{f}_{nh}(x) - \mathrm{E}\,[\hat{f}_{nh}(x)])$$
$$= \frac{A}{\pi} \left(\frac{1}{\sigma^2}\right)^{1+\alpha} (\Gamma(\alpha+1) + o(1)) U_{nh}(x) + O_P(h),$$

*where*

$$U_{nh}(x) = \frac{1}{\sqrt{n}} \sum_{j=1}^n \left( \cos\left(\frac{X_j - x}{h}\right) - \mathrm{E}\left[\cos\left(\frac{X_j - x}{h}\right)\right] \right). \tag{5.3}$$

*Proof.* We have

$$\hat{f}_{nh}(x) = \frac{1}{2\pi} \int_{-\infty}^\infty e^{-itx} \phi_w(ht) \frac{1}{e^{-\sigma^2 t^2/2}} \phi_{emp}(t) dt$$
$$= \frac{1}{2\pi h} \int_{-1}^1 e^{-isx/h} \phi_w(s) e^{\sigma^2 s^2/(2h^2)} \phi_{emp}\left(\frac{s}{h}\right) ds$$
$$= \frac{1}{2\pi nh} \sum_{j=1}^n \int_{-1}^1 e^{is(X_j - x)/h} \phi_w(s) e^{\sigma^2 s^2/(2h^2)} ds$$
$$= \frac{1}{\pi nh} \sum_{j=1}^n \int_0^1 \cos\left(s\left(\frac{X_j - x}{h}\right)\right) \phi_w(s) e^{\sigma^2 s^2/(2h^2)} ds.$$

Notice that
$$\cos\left(s\left(\frac{X_j - x}{h}\right)\right) = \cos\left(\frac{X_j - x}{h}\right) + \left(\cos\left(s\left(\frac{X_j - x}{h}\right)\right) - \cos\left(\frac{X_j - x}{h}\right)\right)$$
$$= \cos\left(\frac{X_j - x}{h}\right) - 2\sin\left(\frac{1}{2}(s+1)\left(\frac{X_j - x}{h}\right)\right)\sin\left(\frac{1}{2}(s-1)\left(\frac{X_j - x}{h}\right)\right)$$
$$= \cos\left(\frac{X_j - x}{h}\right) + R_{n,j}(s),$$

where $R_{n,j}(s)$ is a remainder term satisfying
$$|R_{n,j}| \leq (|x| + |X_j|)\left(\frac{1-s}{h}\right), \quad 0 \leq s \leq 1. \tag{5.4}$$

The bound follows from the inequality $|\sin x| \leq |x|$.

By Lemma 5.1, $\hat{f}_{nh}(x)$ equals
$$\frac{1}{\pi h}\int_0^1 \phi_w(s)e^{\sigma^2 s^2/(2h^2)}ds \frac{1}{n}\sum_{j=1}^n \cos\left(\frac{X_j - x}{h}\right) + \frac{1}{n}\sum_{j=1}^n \tilde{R}_{n,j}$$
$$= \frac{A}{\pi}(\Gamma(\alpha+1) + o(1))\left(\frac{1}{\sigma^2}\right)^{1+\alpha} h^{1+2\alpha}e^{\sigma^2/(2h^2)}\frac{1}{n}\sum_{j=1}^n \cos\left(\frac{X_j - x}{h}\right)$$
$$+ \frac{1}{n}\sum_{j=1}^n \tilde{R}_{n,j},$$

where
$$\tilde{R}_{n,j} = \frac{1}{\pi}\frac{1}{h}\int_0^1 R_{n,j}(s)\phi_w(s)e^{\sigma^2 s^2/(2h^2)}ds.$$

For the remainder we have, by (5.4) and Lemma 5.1,
$$|\tilde{R}_{n,j}| \leq \frac{1}{\pi}(|x| + |X_j|)\frac{1}{h}\int_0^1 \left(\frac{1-s}{h}\right)\phi_w(s)e^{\sigma^2 s^2/(2h^2)}ds$$
$$= \frac{A}{\pi}(|x| + |X_j|)(\Gamma(\alpha+2) + o(1))h^{2+2\alpha}\left(\frac{1}{\sigma^2}\right)^{2+\alpha}e^{\sigma^2/(2h^2)}.$$

Consequently,
$$\mathrm{Var}\left[\tilde{R}_{n,j}\right] \leq \mathrm{E}\left[\tilde{R}_{n,j}^2\right] = O\left(h^{4+4\alpha}e^{\sigma^2/h^2}\right)$$

and
$$\frac{1}{n}\sum_{j=1}^n (\tilde{R}_{n,j} - \mathrm{E}\left[\tilde{R}_{n,j}\right]) = O_P\left(\frac{h^{2+2\alpha}}{\sqrt{n}}e^{\sigma^2/(2h^2)}\right),$$

which follows from Chebyshev's inequality. Finally, we get
$$\frac{\sqrt{n}}{h^{1+2\alpha}e^{\sigma^2/(2h^2)}}(\hat{f}_{nh}(x) - \mathrm{E}\left[\hat{f}_{nh}(x)\right])$$
$$= \frac{A}{\pi}(\Gamma(\alpha+1) + o(1))\left(\frac{1}{\sigma^2}\right)^{1+\alpha} U_{nh}(x) + O_P(h),$$

and this completes the proof of the lemma. □

The next lemma establishes the asymptotic normality.

**Lemma 5.3.** *Assume conditions of Lemma 5.2 and let, for a fixed $x$, $U_{nh}(x)$ be defined by (5.3). Then, as $n \to \infty$ and $h \to 0$,*

$$U_{nh}(x) \xrightarrow{\mathcal{D}} N\left(0, \frac{1}{2}\right).$$

*Proof.* Write

$$Y_j = \frac{X_j - x}{h} \mod 2\pi.$$

For $0 \le y < 2\pi$ we have

$$P(Y_j \le y) = \sum_{k=-\infty}^{\infty} P(2k\pi h + x \le X_j \le 2k\pi h + yh + x)$$

$$= \sum_{k=-\infty}^{\infty} \int_{2k\pi h + x}^{2k\pi h + yh + x} q(u)du = \sum_{k=-\infty}^{\infty} yh\, q(\xi_{k,h})$$

$$= \frac{y}{2\pi} \sum_{k=-\infty}^{\infty} 2\pi h\, q(\xi_{k,h}) \sim \frac{y}{2\pi} \int_{-\infty}^{\infty} q(u)du = \frac{y}{2\pi},$$

where $\xi_{k,h}$ is a point in the interval $[2k\pi h + x, 2k\pi h + yh + x] \subset [2k\pi h + x, 2(k+1)\pi h + x]$. Since $h \to 0$, the last equivalence follows from a Riemann sum approximation of the integral and continuity of the density $q$ of $X$. Consequently, as $h \to 0$, we have $Y_j \xrightarrow{\mathcal{D}} U$, where $U$ is uniformly distributed on the interval $[0, 2\pi]$. Since the cosine is bounded and continuous, it then follows by the dominated convergence theorem that $E\left[|\cos Y_j|^a\right] \to E\left[|\cos U|^a\right]$, for all $a > 0$. Therefore

$$E\left[\cos\left(\frac{X_j - x}{h}\right)\right] \to E\left[\cos U\right] = 0$$

and

$$E\left[\left(\cos\left(\frac{X_j - x}{h}\right)\right)^2\right] \to E\left[(\cos U)^2\right] = \frac{1}{2}.$$

To prove asymptotic normality of $U_{nh}(x)$, first note that it is a normalised sum of i.i.d. random variables. We will verify that the conditions for asymptotic normality in the triangular array scheme of Theorem 7.1.2 in [8] hold (Lyapunov's condition). In our case this reduces to the verification of the fact that

$$\sum_{j=1}^{n} \frac{E\left[|\cos Y_j - E\left[\cos Y_j\right]|^3\right]}{n^{3/2}(\text{Var}[\cos Y_j])^{3/2}} = \frac{E\left[|\cos Y_1 - E\left[\cos Y_1\right]|^3\right]}{n^{1/2}(\text{Var}[\cos Y_1])^{3/2}} \to 0.$$

Now notice that

$$\frac{E\left[|\cos Y_1 - E\left[\cos Y_1\right]|^3\right]}{n^{1/2}(\text{Var}[\cos Y_1])^{3/2}} \sim \frac{E\left[|\cos U|^3\right]}{n^{1/2}(\text{Var}[\cos U])^{3/2}} \to 0,$$

as $n \to \infty$. Consequently, $U_{nh}$ is asymptotically normal,

$$U_{nh}(x) \xrightarrow{\mathcal{D}} N\left(0, \frac{1}{2}\right).$$ □

The following corollary immediately follows from Lemmas 5.2 and 5.3.

**Corollary 5.1.** *Under the conditions of Lemma 5.2 we have that*

$$\frac{\sqrt{n}}{h^{1+2\alpha}e^{\sigma^2/(2h^2)}}\left(\hat{f}_{nh}(x) - \mathrm{E}\left[\hat{f}_{nh}(x)\right]\right) \xrightarrow{\mathcal{D}} N\left(0, \frac{A^2(\Gamma(\alpha+1))^2}{2\pi^2}\left(\frac{1}{\sigma^2}\right)^{2+2\alpha}\right).$$

Now we prove Theorem 2.2.

*Proof of Theorem 2.2.* From (1.4) we have that

$$f_{nh}(x) - \mathrm{E}[f_{nh}(x)] = \frac{1}{1-p}(\hat{f}_{nh}(x) - \mathrm{E}[\hat{f}_{nh}(x)]).$$

Hence the result follows from Corollary 5.1. □

The following lemma gives the order of the variance of $f_{nh}(x)$.

**Lemma 5.4.** *Let Condition 1.2 hold and $f_{nh}(x)$ be defined as in (1.4). Then, as $n \to \infty$ and $h \to 0$,*

$$\mathrm{Var}[f_{nh}(x)] = O\left(\frac{h^{2(1+2\alpha)}e^{\sigma^2/h^2}}{n}\right).$$

*Proof.* We have

$$\mathrm{Var}[f_{nh}(x)] = \frac{1}{4\pi^2(1-p)^2}\frac{1}{nh^2}\mathrm{Var}\left[\int_{-1}^{1}e^{is(X_1-x)/h}\phi_w(s)e^{\sigma^2s^2/(2h^2)}ds\right].$$

Notice that

$$\mathrm{Var}\left[\int_{-1}^{1}e^{is(X_1-x)/h}\phi_w(s)e^{\sigma^2s^2/(2h^2)}ds\right] \leq \left(2\int_{0}^{1}|\phi_w(s)|e^{\sigma^2s^2/(2h^2)}ds\right)^2.$$

Recalling Lemma 5.1, we conclude that

$$\mathrm{Var}[f_{nh}(x)] = O\left(\frac{h^{2(1+2\alpha)}e^{\sigma^2/h^2}}{n}\right).$$ □

Next we deal with consistency of $p_{ng}$ and prove Theorem 2.3.

*Proof of Theorem 2.3.* We have

$$p_{ng} - p = (p_{ng} - \mathrm{E}[p_{ng}]) + (\mathrm{E}[p_{ng}] - p).$$

To prove that this expression converges to zero in probability, it is sufficient to prove that $\mathrm{Var}[p_{ng}] \to 0$ and $\mathrm{E}\,[p_{ng}] - p \to 0$ as $n \to \infty, g \to 0$. We have

$$\mathrm{Var}[p_{ng}] = \pi^2 g^2 \,\mathrm{Var}\left[\frac{1}{2\pi}\int_{-1/g}^{1/g} \frac{\phi_{emp}(t)\phi_k(gt)}{e^{-\sigma^2 t^2/2}}dt\right]$$

$$= \pi^2 g^2 \,\mathrm{Var}[\hat{f}_{ng}(0)].$$

Here it is understood that replacing subindex $h$ by $g$ entails replacement of the smoothing characteristic function $\phi_w$ by $\phi_k$. By Lemma 5.4,

$$\mathrm{Var}[p_{ng}] = O\left(\frac{g^{4+4\alpha}e^{\sigma^2/g^2}}{n}\right). \tag{5.5}$$

This converges to zero due to the condition on $g$. Furthermore,

$$\mathrm{E}\,[p_{nh}] - p = p\left(\frac{1}{2}\int_{-1}^{1} \phi_k(t)dt - 1\right) + (1-p)\frac{g}{2}\int_{-1/g}^{1/g} \phi_f(t)\phi_k(gt)dt. \tag{5.6}$$

The first term here is zero, since $\phi_k$ integrates to 2, while the second term converges to zero, which can be seen upon noticing that $\phi_k$ is bounded, $\phi_f$ is integrable and that this term is bounded by

$$\frac{1-p}{2}g\int_{-\infty}^{\infty}|\phi_f(t)|dt,$$

which converges to zero as $g \to 0$. The last part of the theorem follows from the identity

$$\int_{-1/g}^{1/g} \phi_f(t)\phi_k(gt)dt = g^\gamma \int_{-1/g}^{1/g} t^\gamma \phi_f(t)\frac{\phi_k(gt)}{(gt)^\gamma}dt$$

and Conditions 1.1 and 1.3, because Condition 1.3 implies the existence of a constant $K$, such that $\sup_t |\phi_k(t)t^{-\gamma}| < K$. □

Next we prove asymptotic normality of $p_{ng}$.

*Proof of Theorem 2.4.* The result follows from the definition of $p_{ng}$ and Corollary 5.1, because $p_{ng} = g\pi\hat{f}_{ng}(0)$ essentially is a rescaled version of $\hat{f}_{ng}(0)$. □

Now we are ready to prove Theorem 2.5.

*Proof of Theorem 2.5.* Write

$$\frac{\sqrt{n}}{h^{1+2\alpha}e^{\sigma^2/(2h^2)}}(f^*_{nhg}(x) - \mathrm{E}\,[f^*_{nhg}(x)])$$

$$= \frac{\sqrt{n}}{h^{1+2\alpha}e^{\sigma^2/(2h^2)}}\left(\frac{1}{1-\hat{p}_{ng}}\frac{1}{2\pi}\int_{-\infty}^{\infty} e^{-itx}\frac{\phi_{emp}(t)\phi_w(ht)}{e^{-\sigma^2 t^2/2}}dt\right.$$

$$\left.- \mathrm{E}\left[\frac{1}{1-\hat{p}_{ng}}\frac{1}{2\pi}\int_{-\infty}^{\infty} e^{-itx}\frac{\phi_{emp}(t)\phi_w(ht)}{e^{-\sigma^2 t^2/2}}dt\right]\right)$$

$$-\frac{\sqrt{n}}{h^{1+2\alpha}e^{\sigma^2/(2h^2)}}\frac{1}{h}w\left(\frac{x}{h}\right)\left(\frac{\hat{p}_{ng}}{1-\hat{p}_{ng}} - \mathrm{E}\left[\frac{\hat{p}_{ng}}{1-\hat{p}_{ng}}\right]\right). \tag{5.7}$$

We want to prove that the first term is asymptotically normal, while the second term converges to zero in probability. Application of Slutsky's lemma, see Lemma 2.8 in [41], will then imply that the above expression is asymptotically normal.

First we deal with the second term. We have

$$\frac{\sqrt{n}}{h^{2+2\alpha}e^{\sigma^2/(2h^2)}}w\left(\frac{x}{h}\right)\left(\frac{\hat{p}_{ng}}{1-\hat{p}_{ng}} - \mathrm{E}\left[\frac{\hat{p}_{ng}}{1-\hat{p}_{ng}}\right]\right)$$
$$= w\left(\frac{x}{h}\right)\left(\frac{\sqrt{n}}{h^{2+2\alpha}e^{\sigma^2/(2h^2)}}\frac{g^{2+2\alpha}e^{\sigma^2/(2g^2)}}{\sqrt{n}}\right)$$
$$\times \frac{\sqrt{n}}{g^{2+2\alpha}e^{\sigma^2/(2g^2)}}\left(\frac{\hat{p}_{ng}}{1-\hat{p}_{ng}} - \mathrm{E}\left[\frac{\hat{p}_{ng}}{1-\hat{p}_{ng}}\right]\right). \quad (5.8)$$

Note that Condition 1.4 implies

$$\frac{\sqrt{n}}{h^{2+2\alpha}e^{\sigma^2/(2h^2)}}\frac{g^{2+2\alpha}e^{\sigma^2/(2g^2)}}{\sqrt{n}} = \left(\frac{1+\eta_n}{1+\delta_n}\right)^{1+\alpha}\exp\left(\frac{1}{2}(\delta_n - \eta_n)\log n\right) \to 0.$$

Next we prove that

$$\frac{\sqrt{n}}{g^{2+2\alpha}e^{\sigma^2/(2g^2)}}\left(\frac{\hat{p}_{ng}}{1-\hat{p}_{ng}} - \mathrm{E}\left[\frac{\hat{p}_{ng}}{1-\hat{p}_{ng}}\right]\right) \quad (5.9)$$

is asymptotically normal. Then (5.8) will converge to zero in probability, since convergence to a constant in distribution is equivalent to convergence to the same constant in probability and because $w$ is bounded. We have

$$\frac{\sqrt{n}}{g^{2+2\alpha}e^{\sigma^2/(2g^2)}}(\hat{p}_{ng} - \mathrm{E}\left[\hat{p}_{ng}\right]) \xrightarrow{\mathcal{D}} N\left(0, \frac{C^2(\Gamma(1+\alpha))^2}{2}\left(\frac{1}{\sigma^2}\right)^{2+2\alpha}\right), \quad (5.10)$$

which can be seen as follows:

$$\frac{\sqrt{n}}{g^{2+2\alpha}e^{\sigma^2/(2g^2)}}(\hat{p}_{ng} - \mathrm{E}\left[\hat{p}_{ng}\right]) = \frac{\sqrt{n}}{g^{2+2\alpha}e^{\sigma^2/(2g^2)}}(p_{ng} - \mathrm{E}\left[p_{ng}\right])$$
$$+ \frac{\sqrt{n}}{g^{2+2\alpha}e^{\sigma^2/(2g^2)}}(\hat{p}_{ng} - p_{ng} - \mathrm{E}\left[\hat{p}_{ng} - p_{ng}\right]).$$

Due to Theorem 2.4 the first term here yields the asymptotic normality. We will prove that the second term converges to zero in probability. To this end it is sufficient to prove that

$$\mathrm{Var}\left[\frac{\sqrt{n}}{g^{2+2\alpha}e^{\sigma^2/(2g^2)}}(\hat{p}_{ng} - p_{ng})\right] = \frac{n}{g^{4+4\alpha}e^{\sigma^2/g^2}}\mathrm{Var}\left[\hat{p}_{ng} - p_{ng}\right] \to 0. \quad (5.11)$$

It follows from the definition of $\hat{p}_{ng}$ and Lemma 5.1 that

$$\mathrm{Var}\left[\hat{p}_{ng} - p_{ng}\right] \leq \mathrm{E}\left[(1 - \varepsilon_n - p_{ng})^2 1_{[p_{ng} > 1 - \epsilon_n]}\right]$$
$$\leq (2 + 2K^2 g^{4(1+\alpha)}e^{\sigma^2/g^2})\mathrm{P}(p_{ng} > 1 - \varepsilon_n), \quad (5.12)$$

where $K$ is some constant. This and (5.11) imply that we have to prove

$$n \, \mathrm{P}(p_{ng} > 1 - \varepsilon_n) \to 0. \tag{5.13}$$

Now

$$\mathrm{P}(p_{ng} > 1 - \epsilon_n) = \mathrm{P}(p_{ng} - \mathrm{E}\,[p_{ng}] > 1 - \epsilon_n - \mathrm{E}\,[p_{ng}]). \tag{5.14}$$

Denote $t_n \equiv 1 - \epsilon_n - \mathrm{E}\,[p_{ng}]$ and select $n_0$ so large that for $n \geq n_0$, we have $t_n > 0$. Notice that $t_n \to 1 - p$, which follows from (5.6). The probability in (5.14) is bounded by $\mathrm{P}(|p_{ng} - \mathrm{E}\,[p_{ng}]| > t_n)$. Note that

$$p_{ng} = \sum_{j=1}^{n} \frac{1}{n} \pi k_g \left( \frac{-X_j}{g} \right),$$

with

$$k_g(x) = \frac{1}{2\pi} \int_{-1}^{1} e^{-itx} \phi_k(t) e^{\sigma^2 t^2/(2g^2)} dt.$$

By Lemma 5.1, which is applicable in view of Condition 1.3,

$$\frac{1}{n} \pi k_g \left( \frac{-x}{g} \right) \tag{5.15}$$

is bounded by a constant K, say, times $g^{2+2\alpha} e^{\sigma^2/(2g^2)} n^{-1}$. Hoeffding's inequality, see [29], then yields

$$\mathrm{P}(|p_{ng} - \mathrm{E}\,[p_{ng}]| > t_n) \leq 2 \exp\left(-\frac{2t_n^2}{K^2} \frac{n}{g^{2(2+2\alpha)} e^{\sigma^2/g^2}}\right). \tag{5.16}$$

Since now

$$n \, \mathrm{P}(|p_{ng} - \mathrm{E}\,[p_{ng}]| > t_n) \leq 2n \exp\left(-\frac{2t_n^2}{K^2} \frac{n}{g^{2(2+2\alpha)} e^{\sigma^2/g^2}}\right),$$

it is enough to prove that the term on the right-hand side converges to zero. Taking the logarithm yields

$$\log 2 + \log n - \frac{2t_n^2}{K^2} \frac{n}{g^{2(2+2\alpha)} e^{\sigma^2/g^2}}.$$

This diverges to minus infinity, because the last term dominates $\log n$,

$$\frac{n}{g^{2(2+2\alpha)} e^{\sigma^2/g^2}} \frac{1}{\log n} \to \infty.$$

The latter fact can be seen by taking the logarithm of the left-hand side and using (1.18). We obtain

$$\log n - (2 + 2\alpha) \log g^2 - \frac{\sigma^2}{g^2} - \log \log n$$
$$= -\delta_n \log n + (1 + 2\alpha) \log \log n - (2 + 2\alpha) \log \sigma^2 + (2 + 2\alpha) \log(1 + \delta_n) \to \infty,$$

which follows from (1.18). This in turn proves that (5.10) is asymptotically normal. Since the derivative $(y/(1-y))' \neq 0$, a minor variation of the $\delta$-method then implies that (5.9) is also asymptotically normal (see Theorem 3.8 in [41] for the $\delta$-method). Consequently, the second term in (5.7) converges to zero in probability.

We now consider the first term in (5.7) and want to prove that it is asymptotically normal. Rewrite this term as

$$\frac{\sqrt{n}}{h^{1+2\alpha}e^{\sigma^2/(2h^2)}}\left(\frac{1}{1-p}\frac{1}{2\pi}\int_{-\infty}^{\infty}e^{-itx}\frac{\phi_{emp}(t)\phi_w(ht)}{e^{-\sigma^2t^2/2}}dt\right.$$
$$\left.-\operatorname{E}\left[\frac{1}{1-p}\frac{1}{2\pi}\int_{-\infty}^{\infty}e^{-itx}\frac{\phi_{emp}(t)\phi_w(ht)}{e^{-\sigma^2t^2/2}}dt\right]\right)$$
$$+\frac{\sqrt{n}}{h^{1+2\alpha}e^{\sigma^2/(2h^2)}}\left(\left\{\frac{1}{1-\hat{p}_{ng}}-\frac{1}{1-p}\right\}\frac{1}{2\pi}\int_{-\infty}^{\infty}e^{-itx}\frac{\phi_{emp}(t)\phi_w(ht)}{e^{-\sigma^2t^2/2}}dt\right.$$
$$\left.-\operatorname{E}\left[\left\{\frac{1}{1-\hat{p}_{ng}}-\frac{1}{1-p}\right\}\frac{1}{2\pi}\int_{-\infty}^{\infty}e^{-itx}\frac{\phi_{emp}(t)\phi_w(ht)}{e^{-\sigma^2t^2/2}}dt\right]\right).$$

Thanks to Corollary 5.1 the first summand here is asymptotically normal. We will prove that the second term vanishes in probability. Due to Chebyshev's inequality, it is sufficient to study the behaviour of

$$\frac{\sqrt{n}}{h^{1+2\alpha}e^{\sigma^2/(2h^2)}}\operatorname{E}\left[\left|\frac{\hat{p}_{ng}-p}{(1-\hat{p}_{ng})(1-p)}\frac{1}{2\pi}\int_{-\infty}^{\infty}e^{-itx}\frac{\phi_{emp}(t)\phi_w(ht)}{e^{-\sigma^2t^2/2}}dt\right|\right].$$

By the Cauchy-Schwarz inequality, after taking squares, we can instead consider

$$\frac{n}{h^{2(1+2\alpha)}e^{\sigma^2/h^2}}\operatorname{E}\left[\frac{(\hat{p}_{ng}-p)^2}{(1-\hat{p}_{ng})^2(1-p)^2}\right]$$
$$\times\operatorname{E}\left[\left(\frac{1}{2\pi}\int_{-\infty}^{\infty}e^{-itx}\frac{\phi_{emp}(t)\phi_w(ht)}{e^{-\sigma^2t^2/2}}dt\right)^2\right]. \quad (5.17)$$

Notice that

$$\operatorname{E}\left[\left(\frac{1}{2\pi}\int_{-\infty}^{\infty}e^{-itx}\frac{\phi_{emp}(t)\phi_w(ht)}{e^{-\sigma^2t^2/2}}dt\right)^2\right]$$
$$=\operatorname{Var}\left[\frac{1}{2\pi}\int_{-\infty}^{\infty}e^{-itx}\frac{\phi_{emp}(t)\phi_w(ht)}{e^{-\sigma^2t^2/2}}dt\right]$$
$$+\left(\operatorname{E}\left[\frac{1}{2\pi}\int_{-\infty}^{\infty}e^{-itx}\frac{\phi_{emp}(t)\phi_w(ht)}{e^{-\sigma^2t^2/2}}dt\right]\right)^2.$$

It is easy to see that this expression is of order $h^{-2}$. Indeed, due to Lemma 5.4 the first term in this expression is of order $n^{-1}h^{2(1+2\alpha)}e^{\sigma^2/h^2}$. The fact that this in turn is of lower order than $h^{-2}$ can be seen in the same way as we did with (5.6). For the second term we have

$$\left( \mathrm{E} \left[ \frac{1}{2\pi} \int_{-\infty}^{\infty} e^{-itx} \frac{\phi_{emp}(t)\phi_w(ht)}{e^{-\sigma^2 t^2/2}} dt \right] \right)^2$$
$$= \left( \frac{p}{h} w\left(\frac{x}{h}\right) + (1-p) \frac{1}{2\pi} \int_{-\infty}^{\infty} e^{-itx} \phi_f(t) \phi_w(ht) dt \right)^2,$$

and this is of order $h^{-2}$, because

$$\left| (1-p) \frac{1}{2\pi} \int_{-\infty}^{\infty} e^{-itx} \phi_f(t) \phi_w(ht) dt \right| \leq (1-p) \frac{1}{2\pi} \int_{-\infty}^{\infty} |\phi_f(t)| dt$$

and because $w$ is bounded. Consequently, taking into account (5.17), we have to study

$$\frac{n}{h^{4+4\alpha} e^{\sigma^2/h^2}} \mathrm{E} \left[ \frac{(\hat{p}_{ng} - p)^2}{(1-\hat{p}_{ng})^2 (1-p)^2} \right],$$

or

$$\frac{1}{(1-p)^2} \frac{n}{\epsilon_n^2 h^{2(2+2\alpha)} e^{\sigma^2/h^2}} \mathrm{E}\left[(\hat{p}_{ng} - p)^2\right], \quad (5.18)$$

since $(1-\hat{p}_{ng})^{-2} \leq \epsilon_n^{-2}$. Now

$$\mathrm{E}\left[(\hat{p}_{ng} - p)^2\right] \leq 2\mathrm{E}\left[(\hat{p}_{ng} - p_{ng})^2\right] + 2\mathrm{E}\left[(p_{ng} - p)^2\right].$$

Hence we have to prove that

$$\frac{n}{\epsilon_n^2 h^{2(2+2\alpha)} e^{\sigma^2/h^2}} \mathrm{E}\left[(\hat{p}_{ng} - p_{ng})^2\right] \to 0, \quad (5.19)$$

$$\frac{n}{\epsilon_n^2 h^{2(2+2\alpha)} e^{\sigma^2/h^2}} \mathrm{E}\left[(p_{ng} - p)^2\right] \to 0. \quad (5.20)$$

The first fact essentially follows from the arguments concerning (5.11), since the presence of an additional factor $\epsilon_n^{-2}$ given Condition 1.5 does not affect the arguments used. Indeed, (5.19) will hold true, if we prove that

$$\frac{1}{\epsilon_n^2} \frac{n}{h^{4+4\alpha} e^{\sigma^2/h^2}} g^{4+4\alpha} e^{\sigma^2/g^2} \mathrm{P}(p_{ng} > 1 - \epsilon_n) \to 0.$$

Here we used the definitions of $p_{ng}$ and $\hat{p}_{ng}$ and Lemma 5.1. Now notice that $(g/h)^{4+4\alpha} \to 1$, which follows from Condition 1.4 and that by arguments concerning (5.13) we have $n \mathrm{P}(p_{ng} > 1 - \epsilon_n) \to 0$. Moreover, under Conditions 1.4 and 1.5 we have $\epsilon_n^{-2} e^{\sigma^2/2(1/g^2 - 1/h^2)} \to 0$, which can again be seen by taking the logarithm and verifying that it diverges to minus infinity. This proves (5.19). Next we will prove (5.20). Notice that the latter is in turn implied by

$$\frac{n}{\epsilon_n^2 h^{2(2+2\alpha)} e^{\sigma^2/h^2}} \mathrm{Var}[p_{ng}] + \frac{n}{\epsilon_n^2 h^{2(2+2\alpha)} e^{\sigma^2/h^2}} (\mathrm{E}\left[p_{ng} - p\right])^2 \to 0.$$

The first term here converges to zero by (5.5) and Conditions 1.4 and 1.5. Now we turn to the second term. Taking into account (5.6), we have to study the behaviour of

$$\frac{\sqrt{n}}{\epsilon_n h^{2+2\alpha} e^{\sigma^2/(2h^2)}}(1-p)\int_{-1}^{1}\phi_f\left(\frac{t}{g}\right)\phi_k(t)dt.$$

This can be rewritten as

$$\left(\frac{\sqrt{n}}{\epsilon_n h^{2+2\alpha} e^{\sigma^2/(2h^2)}}\frac{g^{2+2\alpha} e^{\sigma^2/(2g^2)}}{\sqrt{n}}\right)\frac{\sqrt{n}}{g^{2+2\alpha} e^{\sigma^2/(2g^2)}}\int_{-1}^{1}\phi_f\left(\frac{t}{g}\right)\phi_k(t)dt.$$

The factor between the brackets in this expression converges to zero. Therefore it is sufficient to consider

$$\frac{\sqrt{n}}{g^{2+2\alpha} e^{\sigma^2/(2g^2)}}\int_{-1}^{1}\phi_f\left(\frac{t}{g}\right)\phi_k(t)dt.$$

Rewrite this as

$$\frac{\sqrt{n}}{g^{1+2\alpha} e^{\sigma^2/(2g^2)}}g^\gamma \int_{-1/g}^{1/g} t^\gamma \phi_f(t)\frac{\phi_k(gt)}{(gt)^\gamma}dt.$$

Conditions 1.1, 1.3, 1.4 and 1.5 imply that this expression converges to zero, because the integral converges to a constant by the dominated convergence theorem, while

$$\frac{\sqrt{n}}{g^{1+2\alpha} e^{\sigma^2/(2g^2)}}g^\gamma \to 0,$$

which can be seen by taking the logarithm and noticing that it diverges to minus infinity. We obtain

$$\frac{1}{2}\log n + (\gamma - 1 - 2\alpha)\log g - \frac{\sigma^2}{2g^2} = -\frac{\delta_n}{2}\log n + (\gamma - 1 - 2\alpha)\log \sigma$$
$$+ \frac{1+2\alpha-\gamma}{2}\log(1+\delta_n) + \frac{1+2\alpha-\gamma}{2}\log\log n$$
$$\leq (\gamma - 1 - 2\alpha)\log \sigma$$
$$+ \frac{1+2\alpha-\gamma}{2}\log(1+\delta_n) + \frac{1+2\alpha-\gamma}{2}\log\log n \to -\infty, \quad (5.21)$$

which follows from the facts that $\delta_n > 0$ and $1 + 2\alpha - \gamma < 0$. Combination of all these intermediary results completes the proof of the theorem. □

*Proof of Theorem 2.6.* Write

$$\mathrm{E}\left[f^*_{nhg}(x)\right] - f(x) = \{\mathrm{E}\left[f^*_{nhg}(x) - f_{nh}(x)\right]\} + \{\mathrm{E}\left[f_{nh}(x)\right] - f(x)\}. \quad (5.22)$$

Because of (1.7), the second summand in this expression vanishes as $h \to 0$. Next we consider the first summand in (5.22). Using the definitions of $f^*_{nhg}(x)$

and $f_{nh}(x)$, we get

$$\mathrm{E}\left[f_{nhg}^*(x) - f_{nh}(x)\right] = \mathrm{E}\left[\frac{\hat{p}_{ng} - p}{(1 - \hat{p}_{ng})(1 - p)} \frac{1}{2\pi} \int_{-\infty}^{\infty} e^{-itx} \frac{\phi_{emp}(t)\phi_w(ht)}{e^{-\sigma^2 t^2/2}} dt\right]$$
$$- \frac{1}{h} w\left(\frac{x}{h}\right) \mathrm{E}\left[\frac{\hat{p}_{ng} - p}{(1 - \hat{p}_{ng})(1 - p)}\right]. \quad (5.23)$$

By the Cauchy-Schwarz inequality the absolute value of the first summand in this expression is bounded by

$$\left\{\mathrm{E}\left[\frac{(\hat{p}_{ng} - p)^2}{(1 - \hat{p}_{ng})^2(1 - p)^2}\right]\right\}^{1/2}$$
$$\times \left\{\mathrm{E}\left[\left(\frac{1}{2\pi}\int_{-\infty}^{\infty} e^{-itx} \frac{\phi_{emp}(t)\phi_w(ht)}{e^{-\sigma^2 t^2/2}} dt\right)^2\right]\right\}^{1/2}. \quad (5.24)$$

The fact that this term converges to zero follows from (5.17) and subsequent arguments in the proof of Theorem 2.5.

Now we have to study the second summand in (5.23). By the Cauchy-Schwarz inequality and the fact that $(1 - \hat{p}_{ng})^{-2} \leq \epsilon_n^{-2}$, it suffices to consider

$$\frac{1}{(1-p)^2} \frac{1}{\epsilon_n^2} \left(\frac{1}{h} w\left(\frac{x}{h}\right)\right)^2 \mathrm{E}\left[(\hat{p}_{ng} - p)^2\right]$$

instead. The fact that this term converges to zero follows from the arguments concerning (5.18), which were given in the proof of Theorem 2.5. Indeed, the expression above can be rewritten as

$$\left(w\left(\frac{x}{h}\right)\right)^2 \frac{1}{(1-p)^2} \frac{n}{\epsilon_n^2 h^{2(2+2\alpha)} e^{\sigma^2/h^2}} \mathrm{E}\left[(\hat{p}_{ng} - p)^2\right] \frac{h^{2(1+2\alpha)} e^{\sigma^2/h^2}}{n}.$$

Now use arguments concerning (5.19), (5.20) and the facts that $w$ is a bounded function and under Condition 1.4 we have $h^{2(1+2\alpha)} e^{\sigma^2/h^2} n^{-1} \to 0$. This concludes the proof of the first part of the theorem.

Now we prove the second part, an order expansion of the bias $\mathrm{E}\left[f_{nhg}^*(x)\right] - f(x)$ under additional assumptions given in the statement of the theorem. The proof follows the same steps as the proof of the first part of the theorem. Notice that under the condition $f \in \mathcal{H}(\beta, L)$, the second summand in (5.22) is of order $h^\beta$, see Proposition 1.2 in [40]. We have to show then that $h^\beta$ times $\sqrt{n} h^{-1-2\alpha} e^{-\sigma^2/(2h^2)}$ converges to zero. To this end it is sufficient to show that

$$\log\left(h^{\beta-1-2\alpha} \sqrt{n} e^{-\sigma^2/(2h^2)}\right) \to -\infty.$$

This essentially follows from the same argument as (5.21) (with $\gamma$ replaced by $\beta$). Now consider (5.23). Its first term is bounded by (5.24) and we have to show that this term multiplied by $\sqrt{n} h^{-1-2\alpha} e^{-\sigma^2/(2h^2)}$ tends to zero. The arguments from the proof of Theorem 2.6 lead us to (5.18) and hence the desired result. □

*Proof of Theorem 2.7.* The result is a direct consequence of Theorems 2.5 and 2.6. □

## Acknowledgements

The authors would like to thank Chris Klaassen for his careful reading of the draft version of the paper and suggestions that led to improvement of its readability. Comments by an associate editor and two referees are gratefully acknowledged.